\documentclass[aop]{imsart}

 \RequirePackage{amsthm,amsmath,amsfonts,amssymb}
 \RequirePackage[numbers]{natbib}
 \RequirePackage[colorlinks,citecolor=blue,urlcolor=blue]{hyperref}
 \RequirePackage{graphicx}

 \usepackage{tikz} % draw picture
 \usepackage{amsmath}
 \usepackage{mathrsfs} % math notations i.e. \mathscr
 \usepackage{amsfonts} % the font of math notation i.e. \mathbb
 \usepackage{pgfplots}

 \def\ar{\!\!\!&} %definition of \ar

 \startlocaldefs

 \numberwithin{equation}{section}
 \newtheorem{theorem}{Theorem}[section]
 \newtheorem{condition}[theorem]{Condition}
 
 \newtheorem{lemma}[theorem]{Lemma}
 \newtheorem{corollary}[theorem]{Corollary}
 \newtheorem{proposition}[theorem]{Proposition}
 \theoremstyle{remark}
 \newtheorem{remark}[theorem]{Remark}
 
 \endlocaldefs

  \def\beqlb{\begin{eqnarray}}\def\eeqlb{\end{eqnarray}} %Math formulas with numbering
  \def\beqnn{\begin{eqnarray*}}\def\eeqnn{\end{eqnarray*}} %Math formulas without numbering
  \def\ar{&} %definition of \ar

  \def\proof{\noindent{\it Proof.~~}} % definition of Proof.
  \def\qed{\hfill$\Box$\medskip} % definition of the signal of the end of the proof
 
\begin{document}

 \begin{frontmatter}

 \title{Asymptotic Results for Rough Continuous-state Branching Processes}
 \runtitle{Asymptotic Results for Rough CB-processes}

 \begin{aug}
 \author[a]{\fnms{Wei} \snm{Xu}
	\ead[label=e1]{xuwei@math.hu-berlin.de}}
 \address[a]{Department of Mathematics, Humboldt-Universit\"at zu Berlin, Unter den Linden 6, 10099 Berlin, Germany.
 \printead{e1}}
 \runauthor{W. Xu}
% \affiliation{Humboldt-Universit\"at zu Berlin}
 \end{aug}

 \begin{abstract}
  In this paper	 we provide explicit representations of Laplace transforms of extinction time and total progeny of rough continuous-state branching processes introduced in \cite{Xu2021b}.  Also, we show that  their tail distributions are much fatter than those of Feller branching diffusions. 
 \end{abstract}



 \end{frontmatter}

  \section{Introduction}
 \setcounter{equation}{0}

A \textit{rough continuous-state branching process} (rough CB-process) was introduced in \cite{Xu2021b} in the generalization of the classical second Ray-Knight theorem to a spectrally positive stable process $\{\xi(t):t\geq 0\}$ with index $1+\alpha\in(1,2)$ and Laplace exponent 
\beqlb\label{LaplaceExponent}
{\it\Phi}(\lambda)
= b\lambda + c\lambda^{1+\alpha}
= b\lambda + \int_0^\infty \big( e^{-\lambda y}-1+\lambda y \big)\nu_\alpha(dy),
\eeqlb
 where $b\geq 0$, $c>0$ and  the \textit{L\'evy measure}  $\nu_\alpha(dy)$ is given by 
  \beqlb\label{LevyMeasure}
\nu_\alpha(dy):= \frac{c\alpha(\alpha+1)}{\Gamma(1-\alpha)} y^{-\alpha-2}dy. 
\eeqlb
It is described as the unique weak solution to the stochastic Volterra equation 
 \beqlb\label{FractionalSVR}\qquad
X_\zeta(t)= \zeta- b \int_0^t \frac{(t-s)^{\alpha-1}}{c\Gamma(\alpha)}  X_\zeta(s) ds  + \int_0^t \int_0^\infty \int_0^{X_\zeta(s)} \int_{(t-s-y)^+}^{t-s}\frac{r^{\alpha-1}}{c\Gamma(\alpha)} dr
\tilde{N}_\alpha(ds,dy,dz),
\eeqlb 
 where  $\tilde{N}_\alpha(ds,dy,dz)$ is a compensated Poisson random measure on $(0,\infty) \times \mathbb{R}_+^2$ with intensity $ds\nu_\alpha(dy)dz$.
 
 The rough CB-process $X_\zeta$ is locally H\"older-continuous of any order strictly less than $\alpha/2$ and will fall into the trap $0$ in finite time, i.e. its \textit{extinction time} 
 \beqlb
 \tau_{X_\zeta}:= \inf\{ t\geq 0: X_\zeta(t)=0 \}  
 \eeqlb
 is finite almost surely and $X_\zeta(\tau_{X_\zeta}+t)= 0 $ for any $t\geq 0$. Thus its \textit{total progeny} defined by
 \beqlb
 \mathbf{T}_{X_\zeta}:= \int_0^\infty X_\zeta(t)dt = \int_0^{\tau_{X_\zeta}} X_\zeta(t)dt
 \eeqlb
 is also finite almost surely. 
 In this work, we mainly explore the asymptotic behavior of  tail distributions of $\tau_{X_\zeta}$ and $\mathbf{T}_{X_\zeta}$. 
  
 To understand the laws of $\tau_{X_\zeta}$ and $\mathbf{T}_{X_\zeta}$, we first study the characteristic functionals of $X_\zeta$. 
 From Theorem~1.3 in \cite{Xu2021b}, there exist two continuous functions $\{v_t(\lambda): t,\lambda \geq 0\}$ and $\{ V_\mathbf{T}(\lambda):\lambda \geq 0  \}$ such that 
 \beqlb
 \mathbf{E}\big[  \exp\{-\lambda X_
 \zeta(t)\} \big] = \exp\{ -\zeta \cdot v_t(\lambda)  \} 
 \quad \mbox{and}\quad 
   \mathbf{E}\big[  \exp\{-\lambda \mathbf{T}_{X_\zeta}\} \big] = \exp\{ -\zeta\cdot V_\mathbf{T}(\lambda) \} . 
 \eeqlb
 Enlightened by the method developed in the proof of Theorem~1.1 in \cite{Xu2021b}, we first provide explicit representations of $v_t(\lambda)$ and $V_\mathbf{T}(\lambda)$ in terms of functionals of scale function and inverse Laplace exponent of $\xi$. Our main idea is characterizing the limits of characteristic functionals of local time processes of rescaled compounded Poisson processes which converge weakly to the stable process $\xi$. 
 Noting that the function $v_t(\lambda)$ increases to a finite limit $\bar{v}_t$ as $\lambda \to \infty$ and the function $\bar{v}_t$ decreases to $0$ as $t\to \infty$, we have 
 \beqnn
 \mathbf{P}\{ \tau_{X_\zeta} >t \} = \mathbf{P}\{ X_\zeta(t) >0 \} =1- e^{-\zeta\cdot\bar{v}_t} \sim \zeta \cdot\bar{v}_t,
 \eeqnn
 as $t\to\infty$. 
 Applying Karamata-Tauberian theorem  together with  our observation that both $v_t(\lambda)$ and $V_\mathbf{T}(\lambda)$ are regularly varying at infinity, we obtain that both $\tau_\zeta$ and $\mathbf{T}_\zeta$ have heavy-tailed distributions.
 In precise, for both $b=0$  and  $b>0$ we have as $t\to\infty$,
 \beqlb
 \mathbf{P}\{ \tau_{X_\zeta}>t\}  \sim C t^{-\alpha}.
 \eeqlb
 Differently,   as $x\to\infty$ we have
 \beqlb
 \mathbf{P}\{ \mathbf{T}_{X_\zeta}>x  \}  \sim C\cdot x^{-\alpha/(1+\alpha)}\mbox{ when }b=0 \mbox{ and }\mathbf{P}\{ \mathbf{T}_{X_\zeta}>x  \}  \sim C\cdot x^{-\alpha} \mbox{ when }b>0.
 \eeqlb 
 These are in sharp contrast to the asymptotic behavior of Feller branching diffusion defined as the unique solution to 
 \beqlb\label{Eqn.01}
Y_\zeta(t)
\ar=\ar \zeta - \int_0^t b Y_\zeta(s)ds +  \int_0^t \sqrt{2cY_\zeta(s)} dB(s),
\eeqlb
  where $\{B(t):t\geq 0\}$ is a standard Brownian motion. 
  In detail,
  denote by $\tau_{Y_\zeta}$ and $\mathbf{T}_{Y_\zeta}$ its extinction time and total progeny respectively. 
  Theorem~3.4  in \cite[p.17]{Li2019} together with a simply calculation shows that 
  \beqlb
  \mathbf{P}\{ \tau_{Y_\zeta}>t\}  \sim C t^{-1} \mbox{ when }b=0 \mbox{ and } \mathbf{P}\{ \tau_{Y_\zeta}>t\}  \sim Ce^{-bt} \mbox{ when }b>0.
  \eeqlb  
  Moreover, from Corollary~10.9 in \cite[p.280]{Kyprianou2006}  we also have 
    \beqlb
  \mathbf{P}\{ \mathbf{T}_{Y_\zeta}>x  \}  \sim C\cdot x^{-1/2}\mbox{ when }b=0 \mbox{ and }\mathbf{P}\{ \mathbf{T}_{Y_\zeta}>x  \}  \leq C e^{-\frac{b^2}{4c}x} \mbox{ when }b>0.
  \eeqlb

   \smallskip
  \textit{Organization of this paper.} In Section~2, we first give a criticality for a rough CB-process and then provide the asymptotic results for its extinction time and total progeny.
  The proofs are given in Section~3.

  \section{Main results}\label{preliminaries}
  \setcounter{equation}{0}  

% \subsection{Criticality} 
Define the \textit{Mittag-Leffler function} $E_{\alpha,\alpha}$ on $\mathbb{R}_+$ by
\beqnn
E_{\alpha,\alpha}(x):= \sum_{n=0}^\infty \frac{x^n}{\Gamma(\alpha(n+1))}.
\eeqnn
It is locally H\"older-continuous with index $\alpha$; see \cite{HauboldMathaiSaxena2011} for a precise definition of it and a survey of some of its properties, e.g. for any $a\geq 0 $ we have the well-known Laplace transform
\beqlb\label{Mittag-Leffler}
\int_0^\infty  e^{-\lambda x} a x^{\alpha-1} E_{\alpha,\alpha}(-a\cdot x^\alpha) dx = \frac{a}{a+\lambda^{\alpha}} ,  \quad \lambda \geq 0.
\eeqlb 
Let $\{W(x):x\in\mathbb{R}\}$ be the \textit{scale function} of $\xi$,  which is identically zero on $x\in (-\infty,0)$ and characterized on $[0,\infty)$ as a strictly increasing function whose Laplace transform is given by
\beqlb\label{ScaleFunction}
\int_0^\infty e^{-\lambda x}W(x)dx= \frac{1}{\mathit{\Phi}(\lambda )},
\quad \lambda >0.
\eeqlb
The scale function $W$ is infinitely differentiable on $(0,\infty)$ with derivative denoted as $W'$. 
Applying integration by parts to (\ref{ScaleFunction}) and the using (\ref{Mittag-Leffler}), we see that the derivative function $W'$ is of the form:
\beqnn
W'(x)=   c^{-1}x^{\alpha-1}  \cdot E_{\alpha,\alpha}(-b/c\cdot x^\alpha),\quad x\geq 0.
\eeqnn
Specially, when $b=0$ we have $E_{\alpha,\alpha}(0)=\Gamma(\alpha)$ and 
\beqlb\label{ScaleF0}
W'(x)= \frac{x^{\alpha-1}}{c\Gamma(\alpha)},\quad W(x)=\frac{x^{\alpha}}{c\Gamma(1+\alpha)}.
\eeqlb
 When $b>0$, the function $bW'$ turns to be \textit{Mittag-Leffler probability density function} and hence $bW(t)
 \to 1$ as $t\to\infty$.  
  
 From Theorem~1.4 in \cite{Xu2021b}, the stochastic Volterra equation (\ref{FractionalSVR}) is equivalent to 
  \beqlb\label{MainThm.SVE}
 X_\zeta(t)\ar=\ar \zeta (1-bW(t))+ \int_0^t \int_0^\infty \int_0^{X_\zeta(s)} \big( W(t-s) -W(t-s-y)\big) \tilde{N}_\alpha(ds,dy,dz).
 \eeqlb
 Taking expectations on  both sides of the equation above, we can get the following result.
 \begin{lemma}
 	For $t>0$, we have $\mathbf{E}[X_\zeta(t)]=  \zeta (1-bW(t))$, which is identical to $\zeta$ when $b=0$ and decreases to $0$ as $t\to\infty$ when $b>0$.
 	\end{lemma}
 
 We say the rough CB-process $X_\zeta$ is \textit{critical} or \textit{subcritical}   according as $b=0$ or $b>0$, which is consistent to the criticality for Feller branching diffusions; see  Definition~10.4 in \cite[p.277]{Kyprianou2006}. 

\subsection{Extinction time}
We now provide an explicit representation for the Laplace exponent $v_t(\lambda)$ in term of a  functional of $W$ and then explore the tail distribution of extinction time $\tau_{X_\zeta}$. 
%\beqnn
%\mathbf{E}[\exp\{-\lambda X_\zeta(t)\}]= \exp\{ -\zeta v_t(\lambda)\}
%\eeqnn
\begin{theorem}\label{MainThm.01}
	For $\lambda, t\geq 0$, the Laplace exponent $v_t(\lambda)$ is of the form
	\beqnn
	v_t(\lambda)= \frac{c}{\Gamma(1-\alpha)}  \int_0^\infty \frac{ \lambda\cdot \nabla_x W(t)}{ \lambda W(t)+ 1}  \frac{\alpha dx}{x^{\alpha+1}} 
	\eeqnn 
	and as $\lambda\to \infty$ it increases strictly to 
	\beqnn
	\bar{v}_t:= \frac{c}{\Gamma(1-\alpha)} \int_0^\infty \frac{\nabla_x W(t)}{W(t)} \frac{\alpha dx}{x^{\alpha+1}} <\infty.
	\eeqnn
\end{theorem}

\begin{corollary}\label{MainProp.01}
	The function $\bar{v}_t$ decreases to $0$ as $t \to \infty$.  In precise,
	\begin{enumerate}
		\item[(1)]  when $b=0$, we have 
		\beqnn
		\mathbf{P}\{\tau_{X_\zeta}\geq t\} \sim \zeta\cdot \bar{v}_t = \zeta\cdot c\Gamma(1+\alpha) \cdot t^{-\alpha};
		\eeqnn
		
		\item[(2)] when $b>0$,   we have   
		\beqnn
		\mathbf{P}\{\tau_{X_\zeta}\geq t\} \sim \zeta \cdot \bar{v}_t \sim \frac{\zeta \cdot c}{\Gamma(1-\alpha)}  t^{-\alpha}.
		\eeqnn
	\end{enumerate}
	
\end{corollary}

 \begin{remark}\label{MainRemark.01}
   For $\lambda \geq 0$,	from Theorem~3.4 in \cite[p.17]{Li2019} we have when $b=0$,
 	\beqnn
 	\mathbf{E}\big[\exp\{ -\lambda Y_\zeta(t)  \}\big] = \exp\Big\{-\zeta \cdot \frac{\lambda}{1+ c\lambda t} \Big\}
 	\eeqnn
 	and hence as $t \to \infty$, 
 	\beqnn
 	\mathbf{P}\{  \tau_{Y_\zeta} >t  \}=1- \exp\Big\{-\frac{\zeta}{ct}\Big\} \sim \frac{\zeta}{ct} .
 	\eeqnn
 When $b>0$, from Example 2.3 in \cite[p.14]{Li2019} we have 
  	\beqnn
 \mathbf{E}\big[\exp\{ -\lambda Y_\zeta(t)  \}\big] = \exp\Big\{-\zeta \cdot \frac{\lambda e^{-bt}}{1+ c\lambda \cdot b^{-1}(1-e^{-bt}) } \Big\}
 \eeqnn
 and hence as $t \to \infty$, 
 \beqnn
 \mathbf{P}\{  \tau_{Y_\zeta} >t  \}=1- \exp\Big\{-\zeta \cdot \frac{ e^{-bt}}{ c\cdot b^{-1}(1-e^{-bt})}\Big\} \sim \frac{b}{c}e^{-b t}.
 \eeqnn
 	\end{remark}

%  \begin{proposition}
%  	We have $\bar{v}_t \sim \Gamma(1-\alpha)\Gamma(1+\alpha)t^{-\alpha} $ as $t\to 0+$.
%  \end{proposition}

\subsection{Total progeny}
It is obvious that the function $\mathit{\Phi}$ is infinitely differentiable and strictly increasing on $(0,\infty)$ with $\mathit{\Phi}(0)=0$ and $\mathit{\Phi}(\lambda)\to\infty$ as $\lambda\to \infty$.  
This allows us to define its inverse $\mathit{\Psi}$ on $[0,\infty)$ with $\mathit{\Psi}(0)=0$.

\begin{theorem}\label{MainThm.02}
	For $\lambda \geq 0$	we have  $V_{\mathbf{T}}(\lambda)=c \cdot |\Psi(\lambda)|^\alpha $. 
\end{theorem}
%Since $I_{X_\zeta}<\infty$ a.s., thus as $\lambda \to\infty$ we have $\mathbf{E}[\exp\{ -\lambda  I_{X_\zeta} \}]\to 0$  and hence 	The function $V(\lambda)$ increases strictly to infinity as $\lambda \to \infty$.

\begin{corollary}\label{MainProp.02}
	The function $V_\mathbf{T}(\lambda)$ increases strictly to infinity as $\lambda \to \infty$. In precise,
	\begin{enumerate}
		\item[(1)]  when $b=0$, we have $V_\mathbf{T}(\lambda)= c^{1/(1+\alpha)}  \cdot  \lambda^{\alpha/(1+\alpha)} $ and hence as $x\to\infty$, 
		\beqnn
		\mathbf{P}\{    \mathbf{T}_{X_\zeta}>x \} \sim \frac{\zeta\cdot c^{1/(1+\alpha)}}{\Gamma(1/(1+\alpha))}   \cdot  x^{-\alpha/(1+\alpha)};
		\eeqnn 
		
		\item[(2)] 	when $b>0$, we have $V_\mathbf{T}(\lambda)\sim \zeta\cdot c \cdot (\lambda /b )^\alpha$ as $\lambda \to 0+$ and hence as $x\to\infty$,
		\beqnn
		\mathbf{P}\{    \mathbf{T}_{X_\zeta}>x \} \sim  \frac{\zeta\cdot c}{b^\alpha\Gamma(1-\alpha)} \cdot x^{-\alpha}. 
		\eeqnn 
		
	\end{enumerate}
	
\end{corollary}

\begin{remark}\label{MainRemark.02}
From Corollary~10.9 in \cite[p.280]{Kyprianou2006},  for $\lambda \geq -\frac{b^2}{4c}$ we have 
	\beqlb
	\mathbf{E}\big[\exp\{ -\lambda \mathbf{T}_{Y_\zeta}  \}\big]
	=\exp\Big\{ - \zeta\cdot  \frac{\sqrt{b^2+4c\lambda}-b}{2c} \Big\}.
	\eeqlb
	When $b=0$, we have $	\mathbf{E}\big[\exp\{ -\lambda \mathbf{T}_{Y_\zeta}  \}\big] \sim 1- \zeta \cdot\sqrt{t/c} $ as $\lambda \to 0+$. Applying  Karamata-Tauberian theorem; see Theorem~8.1.6 in \cite[p.333]{BinghamGoldieTeugels1987}, we have as $x\to\infty$, 
	\beqnn
	\mathbf{P}\{  \mathbf{T}_{Y_\zeta}  >x \} \sim \zeta \cdot \frac{x^{-1/2}}{\sqrt{c}\Gamma(1/2)}. 
	\eeqnn
  When $b>0$, by Chebyshev's inequality we have 
  \beqnn
  	\mathbf{P}\{  \mathbf{T}_{Y_\zeta}  >x \} \ar=\ar 	\mathbf{P}\{ e^{  \frac{b^2}{4c}\mathbf{T}_{Y_\zeta}}  >e^{\frac{b^2}{4c}x} \} \leq e^{-\frac{b^2}{4c}x} \mathbf{E}[ e^{  \frac{b^2}{4c}\mathbf{T}_{Y_\zeta}}  ] \leq \exp\Big\{\frac{b\zeta}{2c}-\frac{b^2}{4c}x\Big\} . 
  \eeqnn

	\end{remark}

  \section{Proofs}\label{Proofs}
\setcounter{equation}{0}

 \subsection{Compound Poisson Processes}\label{CompoundPoissonP}

Let $\Lambda$ be a Pareto II distribution on $\mathbb{R}_+$ with location $0$ and shape $\alpha+1$, i.e.
\beqnn
\Lambda (dx)= (\alpha+1)  (1+x )^{-\alpha-2}dx.
\eeqnn
For $n\geq 1$, let $\{\xi^{(n)}(t):t\geq 0\}$ be a compound Poisson process with drift $-1$, arrival rate $\gamma_n>0$, jump size distribution $\Lambda$ and initial state $\xi^{(n)}(0)$ distributed as 
\beqlb\label{SBLambda}
\Lambda^*(dx):= \alpha \bar\Lambda(x) dx= \alpha (1+x)^{-1-\alpha}dx.
\eeqlb
Here $\bar\Lambda(x):=\Lambda[x,\infty)$ is the tail distribution of $\Lambda$. 
It is a spectrally positive L\'evy process with Laplace exponent %${\it\Phi}^{(n)}$ given by
\beqnn
{\it\Phi}^{(n)}(\lambda):= \lambda + \int_0^\infty (e^{-\lambda x}-1) \gamma_n \Lambda (dx),\quad \lambda \geq 0.
\eeqnn
We are interested in the case $\mathbf{E}[\xi^{(n)}(1)]= \gamma_n/\alpha-1\leq 0$ in which the  function ${\it\Phi}^{(n)}$ increases strictly to infinity. 

Let $c_0= (c/\Gamma(1-\alpha))^{\frac{1}{1+\alpha}}$. 
We define the \textit{first passage time} of $\xi^{(n)}$ in $(-\infty,0]$ by
\beqlb
\tau_0^{(n)}:=\inf\{t\geq 0: \xi^{(n)}_t\leq 0\}. 
\eeqlb
Let $\{L_{\xi^{(n)}}(y,t):y\in\mathbb{R}, t\geq 0 \}$ be the local time process of $\xi^{(n)}$ with
\beqnn
L_{\xi^{(n)}}(y,t):= \# \{s\in(0,t]: \xi^{(n)}(s)=y \}.
\eeqnn
We have $L_{\xi^{(n)}}(y,\infty)=\infty$ a.s. if and only if $\mathbf{E}[\xi^{(n)}(1)]=0$; equivalently, if and only if $\gamma_n=\alpha$. 
Let $\{Z^{(n)}_k(t):t\geq 0\}_{k\geq 1}$ be a sequence of i.i.d. copies of local time process $\{L_{\xi^{(n)}}(t,\tau_0^{(n)}):t\geq 0\}$.
We now consider under the following condition, the weak convergences of rescaled processes $\{\xi_0^{(n)}(t):t\geq 0 \}$ and $\{ S_{\zeta}^{(n)}(t):t\geq 0   \}$ with 
\beqlb\label{ScaledLTP}
\xi_0^{(n)}(t):=\frac{c_0}{n}\xi^{(n)}(n^{1+\alpha}t) 
\quad\mbox{and}\quad 
S_{\zeta}^{(n)}(t):= \frac{1}{c_0n^{\alpha}}\sum_{k=1}^{\lfloor c_0n^\alpha \zeta \rfloor} Z^{(n)}_k(nt/c_0) . 
\eeqlb

\begin{condition}\label{Main.Condition}
	Assume that $n^\alpha(1-\gamma_n/\alpha)\to b/c_0 $ as $n\to\infty$.
\end{condition}

  It is easy to see that the L\'evy process $\xi_0^{(n)}$ has Laplace exponent $n^{1+\alpha}{\it\Phi}^{(n)}(c_0\lambda/n)$ for $\lambda \geq 0$.
Under Condition~\ref{Main.Condition}, a routine computation shows that $n^{1+\alpha}{\it\Phi}^{(n)}(c_0\lambda/n)\to {\it \Phi}(\lambda) $ as $n\to\infty$.
By Corollary~4.3 in \cite[p.440]{JacodShiryaev2003}, we have $\xi_0^{(n)}\to \xi$  weakly in $\mathbf{D}([0,\infty);\mathbb{R})$.
Denote by $\{W^{(n)}(x):x\in\mathbb{R}\}$ the scale of $\xi^{(n)}$. We have $W^{(n)}(0)=1$ and  as $n\to\infty$, 
\beqlb\label{ConvergenceScaleF}
n^{-\alpha}W^{(n)}(nt/c_0)\to c_0 W(t),\quad t\in\mathbb{R}. 
\eeqlb
The next lemma follows directly from Lemma~2.3 and Proof of Theorem~1.1 in \cite{Xu2021b}.

\begin{lemma}\label{WeakConvergence}
 Under Condition~\ref{Main.Condition}, for any $\zeta >0$ we have 	$S_{\zeta}^{(n)}$ converges weakly to the unique weak solution $X_\zeta$ to (\ref{FractionalSVR}) in $\mathbf{D}([0,\infty);\mathbb{R}_+)$ as $n\to\infty$. 
	\end{lemma}

\subsection{Proofs for Theorem~\ref{MainThm.01} and  Corollary~\ref{MainProp.01}} 

From Lemma~\ref{WeakConvergence},  it is necessary to consider the  distribution of $Z_1^{(n)}$ at first. 
For $x>0$, denote by $\mathbf{P}_x$ and $\mathbf{E}_x$ the probability law and expectation of $\xi^{(n)}$ conditioned on $\xi^{(n)}(0)=x$. 

\begin{proposition}\label{Prop.3.1}
	For $n\geq 1$ and $x,t> 0$ with $x\neq t$, we have 
	\beqlb
	\mathbf{E}_t \big[\exp\{-\lambda Z_1^{(n)}(t)
	\}  \big] \ar=\ar   \frac{\frac{W^{(n)}(0)}{W^{(n)}(t)}}{1-e^{-\lambda}(1-\frac{W^{(n)}(0)}{W^{(n)}(t)})} , \label{eqn3.08} \\
	\mathbf{E}_x \big[\exp\{-\lambda Z_1^{(n)}(t)
	\}  \big] 
	\ar=\ar 1- \frac{\nabla_x W^{(n)}(t)}{W^{(n)}(t)}   \frac{1-e^{-\lambda}}{1-e^{-\lambda}(1-\frac{W^{(n)}(0)}{W^{(n)}(t)})}.
	\label{eqn3.09}
	\eeqlb
	
	\end{proposition}
\proof %Noting that $\xi^{(n)}$ is either recurrent or drifting to $-\infty$. If $x>t$, we have $W^{(n)}(t-x)=0$ and  $\mathbf{P}_x\{  Z_1^{(n)}(t=1 \}=1$. Thus (????) holds. We now prove (???) with $x\leq t$. 
For $t>0$, let $\tau_t^+:=\inf\{ s\geq 0: \xi^{(n)}(s)>t \}$ and $T_t:= \inf\{   s\geq 0: \xi^{(n)}(s)=t \}$. Since $\xi^{(n)}$ has neither diffusion nor nonnegative jumps, we have $\tau_t^+ < T_t$ a.s. %Moreover, conditioned on $\tau_t^+<\infty$ we also have $T_t< \infty $. For simplicity, we write $\tau_0^-=\tau_0(\xi^{(n)})$. 
Thus 
\beqlb\label{eqn3.10}
\mathbf{E}_x \big[\exp\{-\lambda Z_1^{(n)}(t)
\}  \big] 
\ar=\ar \mathbf{P}_x\{  \tau_0^{(n)}<  \tau_t^+\} + \mathbf{E}_x \big[\exp\{-\lambda Z_1^{(n)}(t)
\};   \tau_0^{(n)}> \tau_t^+ \big] 
\eeqlb
By the strong Markov property of $\xi^{(n)}$, we have 
\beqlb\label{eqn3.11}
\mathbf{E}_x \big[\exp\{-\lambda Z_1^{(n)}(t)
\};   \tau_0^{(n)}> \tau_t^+ \big] = e^{-\lambda} \mathbf{E}_t \big[\exp\{-\lambda Z_1^{(n)}(t)
\}\big] \mathbf{P}_x\{  \tau_t^+<\tau_0^{(n)} \}.
\eeqlb
Solving (\ref{eqn3.10})-(\ref{eqn3.11}) with $x=t$, we have 
\beqnn
\mathbf{E}_t \big[\exp\{-\lambda Z_1^{(n)}(t)
\}  \big] \ar=\ar \frac{\mathbf{P}_t\{  \tau_0^-<  \tau_t^+\} }{1-e^{-\lambda}   \mathbf{P}_t\{  \tau_t^+<\tau_0^- \}}.
%= \frac{\frac{W^{(n)}(0)}{W^{(n)}(t)}}{1-e^{-\lambda}(1-\frac{W^{(n)}(0)}{W^{(n)}(t)})} .
\eeqnn
Plugging this back into (\ref{eqn3.10})-(\ref{eqn3.11}), for $x\neq t$ we have 
\beqnn
\mathbf{E}_x \big[\exp\{-\lambda Z_1^{(n)}(t)
\}  \big] 
\ar=\ar \mathbf{P}_x\{  \tau_0^-<  \tau_t^+\} + \mathbf{P}_x\{  \tau_t^+<\tau_0^- \} \frac{e^{-\lambda} \mathbf{P}_t\{  \tau_0^-<  \tau_t^+\} }{1-e^{-\lambda}   \mathbf{P}_t\{  \tau_t^+<\tau_0^- \}} .
\eeqnn
From Theorem~8.1(iii) in \cite[p.215]{Kyprianou2006}, we have $\mathbf{P}_t\{  \tau_0^-<  \tau_t^+\}=W^{(n)}(0)/W^{(n)}(t)$ and $\mathbf{P}_x\{  \tau_0^-<  \tau_t^+\}=W^{(n)}(t-x)/W^{(n)}(t)$. 
With a simple calculation we can get the desired result.
\qed
 %(\ref{eqn3.08}) follows. 
%Plugging this back into (\ref{eqn3.10})-(\ref{eqn3.10}), for $x\neq t$ we have 
%\beqnn
%\mathbf{E}_x \big[\exp\{-\lambda Z_1^{(n)}(t)
%\}  \big] 
%\ar=\ar \mathbf{P}_x\{  \tau_0^-<  \tau_t^+\} + e^{-\lambda} \mathbf{E}_t \big[\exp\{-\lambda Z_1^{(n)}(t)
%\}\big] \mathbf{P}_x\{  \tau_t^+<\tau_0^- \}\cr
%\ar=\ar \frac{W^{(n)}(t-x)}{W^{(n)}(t)} + \Big(  1-\frac{W^{(n)}(t-x)}{W^{(n)}(t)}  \Big) \frac{e^{-\lambda}\frac{W^{(n)}(0)}{W^{(n)}(t)}}{1-e^{-\lambda}(1-\frac{W^{(n)}(0)}{W^{(n)}(t)})} \cr
%\ar=\ar 1+ \frac{W^{(n)}(t-x)}{W^{(n)}(t)} -1 +  \frac{\nabla_x W^{(n)}(t)}{W^{(n)}(t)} \frac{e^{-\lambda }\frac{W^{(n)}(0)}{W^{(n)}(t)}}{1-e^{-\lambda}(1- \frac{W^{(n)}(0)}{W^{(n)}(t)})} \cr
%\ar=\ar 1- \frac{\nabla_x W^{(n)}(t)}{W^{(n)}(t)}   \Big( 1- \frac{e^{-\lambda }\frac{W^{(n)}(0)}{W^{(n)}(t)}}{1-e^{-\lambda}(1- \frac{W^{(n)}(0)}{W^{(n)}(t)})} \Big)
%\eeqnn
%Applying Theorem~8.1(iii) in \cite[p.215]{Kyprianou2006} again, we have and 
%\qed

%\begin{corollary}
%	For $t\geq 0$, we have 	\beqnn
%	\mathbf{E}\big[\exp\{-\lambda Z_1^{(n)}(t)\}  \big]  
%	\ar=\ar 1- \frac{1-e^{-\lambda}}{1-e^{-\lambda}(1- \frac{W^{(n)}(0)}{W^{(n)}(t)})} \int_0^\infty \frac{\nabla_x W^{(n)}(t)}{W^{(n)}(t)} \Lambda^*(dx)  %\frac{\alpha dx}{(1+x)^{\alpha+1}} 
%	\eeqnn
%\end{corollary}

\textit{Proof for Theorem~\ref{MainThm.01}.} From Proposition~\ref{Prop.3.1}, we have 
	\beqnn
\mathbf{E}\big[\exp\{-\lambda Z_1^{(n)}(t)\}  \big]  
\ar=\ar 1- \frac{1-e^{-\lambda}}{1-e^{-\lambda}(1- \frac{W^{(n)}(0)}{W^{(n)}(t)})} \int_0^\infty \frac{\nabla_x W^{(n)}(t)}{W^{(n)}(t)} \Lambda^*(dx)  %\frac{\alpha dx}{(1+x)^{\alpha+1}} 
\eeqnn
By the mutual independence among $\{Z^{(n)}_k\}_{n\geq 1}$, 
\beqnn
\mathbf{E}_x \big[\exp\{-\lambda S_{\zeta}^{(n)}(t)\}  \big] 
\ar=\ar \Big(  \mathbf{E} \Big[\exp\Big\{-\frac{\lambda}{c_0n^{\alpha}} Z_1^{(n)}(nt/c_0)\Big\}  \Big]      \Big)^{\lfloor n^\alpha c_0\zeta \rfloor}\cr
\ar=\ar \exp\Big\{-\frac{\lfloor n^\alpha c_0\zeta \rfloor}{c_0n^\alpha} v^{(n)}_t(\lambda)\Big\},
%\ar=\ar \Big[ 1- \frac{1-e^{-\frac{\lambda}{c_0n^{\alpha}}}}{1-e^{-\frac{\lambda}{c_0n^{\alpha}}}(1- \frac{W^{(n)}(0)}{W^{(n)}(nt/c_0)})}  \int_0^\infty \frac{\nabla_x W^{(n)}(nt/c_0)}{W^{(n)}(nt/c_0)} \frac{\alpha dx}{(1+x)^{\alpha+1}} \Big]^{\lfloor n^\alpha c_0\zeta \rfloor}.
\eeqnn
where
\beqnn
v^{(n)}_t(\lambda)
\ar=\ar -c_0n^{\alpha}\cdot \log\Big( 1- \frac{1-e^{-\frac{\lambda}{c_0n^{\alpha}}}}{1-e^{-\frac{\lambda}{c_0n^{\alpha}}}(1- \frac{W^{(n)}(0)}{W^{(n)}(nt/c_0)})}  \int_0^\infty \frac{\nabla_x W^{(n)}(nt/c_0)}{W^{(n)}(nt/c_0)} \Lambda^*(dx) \Big). 
%\ar\sim \ar n^{\alpha}   \frac{1-e^{-n^{-\alpha}\lambda}}{1-e^{-n^{-\alpha}\lambda}(1- \frac{W^{(n)}(0)}{W^{(n)}(nt)})}   \int_0^\infty \frac{\nabla_x W^{(n)}(nt) }{W^{(n)}(nt)} \frac{\alpha dx}{(1+x)^{\alpha+1}} 
\eeqnn
From Lemma~\ref{WeakConvergence} we have $v^{(n)}_t(\lambda) \to v_t(\lambda)$ as $n\to\infty$. 
Noting that 
\beqnn
v^{(n)}_t(\lambda)\sim \frac{c_0n^{\alpha}(1-e^{-\frac{\lambda}{c_0n^{\alpha}}})}{c_0n^\alpha\big[1-e^{-\frac{\lambda}{c_0n^{\alpha}}}\big(1- \frac{W^{(n)}(0)}{W^{(n)}(nt/c_0)}\big)\big]}  \int_0^\infty \frac{\nabla_x W^{(n)}(nt/c_0)\Lambda^*(dx)}{n^{-\alpha}W^{(n)}(nt/c_0)/c_0}  . 
\eeqnn
From (\ref{ConvergenceScaleF}) and the fact that $W^{(n)}(0)=1$, we have as $n\to\infty$
\beqnn
 c_0n^\alpha \Big[  1-e^{-\frac{\lambda}{c_0n^{\alpha}}}\Big(1- \frac{W^{(n)}(0)}{W^{(n)}(nt/c_0)}\Big) \Big] 
 \to \lambda +1/W(t).
\eeqnn
and 
\beqnn
 \int_0^\infty \nabla_x W^{(n)}(nt/c_0)  \Lambda^*(dx) 
 \ar=\ar \int_0^\infty \nabla_{nx} W^{(n)}(nt/c_0)  \Lambda^*(n\cdot dx) \to c_0^{1+\alpha} \int_0^\infty \nabla_{x}W(t) \frac{\alpha dx}{x^{1+\alpha}}.
\eeqnn
Putting all estimates above together, we can immediately get the desire result. 
\qed
%
%\beqnn
%c_0n^\alpha(1-e^{-\frac{\lambda}{c_0n^{\alpha}}}) \sim \lambda
%\eeqnn 
%and 
%\beqnn
%n^{-\alpha} W^{(n)}(nt/c_0) \to c_0W(t)  ,\quad W^{(n)}(0)=1
%\eeqnn
%and 
%\beqnn
%n^\alpha(1-e^{-\frac{\lambda}{c_0n^{\alpha}}}) + e^{-\frac{\lambda}{c_0n^{\alpha}}}\cdot n^\alpha \frac{W^{(n)}(0)}{W^{(n)}(nt/c_0)}\sim \lambda/c_0 + \frac{1}{c_0W(t)}.
%\eeqnn
%
%
%Moreover,
%\beqnn
%\int_0^\infty (W^{(n)}(nt/c_0)-W^{(n)}(nt/c_0-x)) \frac{\alpha dx}{(1+x)^{\alpha+1}} 
%\ar=\ar   \int_0^\infty (W^{(n)}(nt/c_0)-W^{(n)}(n(t/c_0-x))) \frac{n\alpha dx}{(1+nx)^{\alpha+1}}\cr
%\ar\sim\ar   \int_0^\infty  n^{-\alpha}(W^{(n)}(nt/c_0)- W^{(n)}(n(t/c_0-x)))\frac{\alpha dx}{x^{\alpha+1}}\cr
%\ar\sim\ar  \int_0^\infty c_0(W(t)- W(t-c_0x))\frac{\alpha dx}{x^{\alpha+1}}\cr
%\eeqnn
%
%Putting them together, we have 
%\beqnn
%v^{(n)}_t(\lambda)\to v_t(\lambda)
%\ar:=\ar   \frac{ \lambda}{ \lambda W(t)+ 1}   \int_0^\infty (W(t)-W(t-c_0x)) \frac{\alpha dx}{x^{\alpha+1}} \cr
%\ar=\ar  \frac{c_0^\alpha \lambda}{ \lambda W(t)+ 1}   \int_0^\infty (W(t)-W(t-x)) \frac{\alpha dx}{x^{\alpha+1}} 
%\eeqnn
 
 \textit{Proof for Corollary~\ref{MainProp.01}.}  
 Since $X_\zeta$ falls into the trap $0$ in finite time, we have $\bar{v}_t\to0 $ as $t\to\infty$.   
 Moreover, we have $\mathbf{P}\{  \tau_{X_{\zeta}}>t \} = 1-\exp\{-\zeta \cdot \bar{v}_t \} \sim \zeta \cdot \bar{v}_t$ as $t\to\infty$.  
Changing the order of integration, we have 
\beqlb\label{eqn3.15}
\int_0^\infty \nabla_xW(t) \frac{\alpha dx}{x^{\alpha+1}}
%\ar=\ar   \int_0^\infty \int_{t-x}^tW'(s)ds \frac{\alpha dx}{x^{\alpha+1}}\cr
\ar=\ar \int_0^\infty \int_0^xW'(t-s)ds \frac{\alpha dx}{x^{\alpha+1}}\cr
\ar=\ar \int_0^\infty W'(t-x)dx  \int_x^\infty \frac{\alpha ds}{s^{\alpha+1}}
= \int_0^t W'(t-x)x^{-\alpha} dx . 
\eeqlb
When $b=0$, from (\ref{ScaleF0}) we have 
\beqnn
\int_0^\infty \nabla_xW(t) \frac{\alpha dx}{x^{\alpha+1}} = \int_0^t \frac{(t-x)^{\alpha-1}x^{-\alpha} dx  }{c\Gamma(\alpha)}  = \frac{\Gamma(1-\alpha)}{c} 
\eeqnn
 and 
 \beqnn
 \mathbf{P}\{  \tau_{X_{\zeta}}>t \}  \sim \zeta /W(t)= W(t)= \zeta\cdot c\Gamma(1+\alpha)\cdot t^{-\alpha}.
 \eeqnn
% \proof we have 
% \beqnn
% \mathbf{P}\{  X_{\zeta}(t)>0 \} 
% \ar=\ar 1- \lim_{\lambda\to\infty}\exp\Big\{ -\zeta \frac{ c_0^\alpha\lambda}{ \lambda W(t)+ 1}   \int_0^\infty (W(t)-W(t-x)) \frac{\alpha dx}{x^{\alpha+1}}  \Big\}\cr
% \ar=\ar 1- \exp\Big\{ -\zeta \frac{c_0^\alpha}{W(t)}   \int_0^\infty (W(t)-W(t-x)) \frac{\alpha dx}{x^{\alpha+1}}  \Big\}\cr
% \ar\sim\ar \zeta \frac{c_0^\alpha}{W(t)}   \int_0^\infty (W(t)-W(t-x)) \frac{\alpha dx}{x^{\alpha+1}}
% \eeqnn 
% {\bf Step 01.} When $b=0$, we have 
% \beqnn
% W(t)= \frac{t^\alpha}{c\Gamma(1+\alpha)} , \quad W'(t) =\frac{t^{\alpha-1}}{c\Gamma(\alpha)} 
% \eeqnn
% and 
% \beqnn
% \int_0^\infty (W(t)-W(t-x)) \frac{\alpha dx}{x^{\alpha+1}}
% \ar=\ar   \int_0^\infty \int_{t-x}^tW'(s)ds \frac{\alpha dx}{x^{\alpha+1}}\cr
% \ar=\ar \int_0^\infty \int_0^xW'(t-s)ds \frac{\alpha dx}{x^{\alpha+1}}\cr
% \ar=\ar \int_0^\infty W'(t-x)dx  \int_x^\infty \frac{\alpha ds}{s^{\alpha+1}}\cr
% \ar=\ar  \int_0^\infty W'(t-x)x^{-\alpha} dx  \cr
% \ar=\ar  \int_0^t \frac{(t-x)^{\alpha-1}x^{-\alpha} dx  }{c\Gamma(\alpha)}  = \frac{\Gamma(1-\alpha)}{c} 
% \eeqnn
% and 
% \beqnn
% \mathbf{P}\{  X_{\zeta}(t)>0 \}  \sim \zeta\cdot c_0^\alpha \cdot  \Gamma(1-\alpha)\Gamma(1+\alpha) \cdot t^{-\alpha}. 
% \eeqnn
 When $b>0$, from (2.4) in \cite{Xu2021b} we have as $t\to\infty$, 
 \beqnn
 W'(t)  
 \sim \frac{c\alpha\cdot t^{-\alpha-1} }{b^2\Gamma(1-\alpha)}=o(t^{-\alpha}).
 \eeqnn
 Applying Proposition 3.11 in \cite{Xu2021a} to (\ref{eqn3.15}), we have as $t\to\infty$,   
 \beqnn
\int_0^\infty \nabla_xW(t) \frac{\alpha dx}{x^{\alpha+1}} 
 \sim W(t) \cdot t^{-\alpha} 
 \eeqnn
 and hence $ \mathbf{P}\{  \tau_{X_{\zeta}}>t \}  \sim \zeta \cdot  c_0 \cdot t^{-\alpha}$.  
 \qed

\subsection{Proofs for Theorem~\ref{MainThm.02} and  Corollary~\ref{MainProp.02}} 

 From Lemma~\ref{WeakConvergence}, we have 
 \beqnn
 \int_0^\infty S_\zeta^{(n)}(t)dt  \to \mathbf{T}_{X_\zeta}
 \eeqnn
 in distribution. From (\ref{ScaledLTP}),
 \beqnn
 \int_0^\infty S_\zeta^{(n)}(t)dt = \sum_{k=1}^{\lfloor c_0n^\alpha \zeta \rfloor}\frac{1}{n^{1+\alpha}}\int_0^\infty  Z^{(n)}_k(t) dt.
 \eeqnn
 Applying the occupation density formula, we immediately have
 \beqnn
 \int_0^\infty Z^{(n)}_1(t)dt= 
 \int_0^\infty L_{\xi^{(n)}}(t,\tau_0^{(n)})dt = \tau_0^{(n)}.
 \eeqnn
 The next proposition can be found in \cite[p.212]{Kyprianou2006}.
 \begin{proposition}\label{Prop3.3}
 Let $\mathit{\Psi}^{(n)}$ be the inverse of $\mathit{\Phi}^{(n)}$. 	For $\lambda \geq 0$, we have 
 	\beqnn
 	\mathbf{E}_x\big[\exp\{-\lambda \tau_0^{(n)}\} \big] = \exp\{ -x  \mathit{\Psi}^{(n)}(\lambda) \}. 
 	\eeqnn
 	
 	\end{proposition}
 
 \textit{Proof for Theorem~\ref{MainThm.02}.}
 By the mutual independence among $\{Z^{(n)}_k\}_{k\geq 1}$ and Proposition~\ref{Prop3.3}, for $\lambda \geq 0$ we have 
 \beqnn
 \mathbf{E}\Big[ \exp\Big\{ -\lambda  \int_0^\infty S_\zeta^{(n)}(t)dt \Big\} \Big]\ar=\ar \Big(   \mathbf{E}\Big[ \exp\Big\{ - \frac{\lambda}{n^{1+\alpha}}\int_0^\infty  Z^{(n)}_k(t) dt\Big\} \Big]   \Big)^{\lfloor c_0n^\alpha \zeta \rfloor}\cr
 \ar=\ar \exp\Big\{-\frac{\lfloor c_0n^\alpha \zeta \rfloor}{c_0n^\alpha} V^{(n)}_\mathbf{T}(\lambda) \Big\},
 \eeqnn
 where 
 \beqnn
 V^{(n)}_\mathbf{T}(\lambda)  \ar=\ar c_0n^\alpha \log \Big(1- \int_0^\infty  \big(1-e^{-x \mathit{\Psi}^{(n)}(\lambda/n^{1+\alpha})} \big)\Lambda^*(dx) \Big).
 \eeqnn
 It is obvious that $V^{(n)}_\mathbf{T}(\lambda) \to V_\mathbf{T}(\lambda)$ as $n\to\infty$. 
 Since $n^{1+\alpha}{\it\Phi}^{(n)}(c_0\lambda/n)\to {\it \Phi}(\lambda) $ as $n\to\infty$, we have $ n\mathit{\Psi}^{(n)}(\lambda/n^{1+\alpha})\to c_0\mathit{\Psi}(\lambda )$. Thus by the dominated convergence theorem,
 \beqnn
 V^{(n)}_\mathbf{T}(\lambda) 
 \ar\sim\ar c_0n^\alpha \int_0^\infty  \big(1-e^{-x \mathit{\Psi}^{(n)}(\lambda/n^{1+\alpha})} \big)\Lambda^*(dx)\cr
 \ar=\ar c_0\int_0^\infty \Big(1- e^{-x\cdot n \mathit{\Psi}^{(n)}(\lambda/n^{1+\alpha})} \Big) \frac{\alpha   n^{\alpha+1}   dx}{(1+nx)^{1+\alpha}}\cr
 \ar\sim \ar  c_0\int_0^\infty \Big(1- e^{-x\cdot c_0\mathit{\Psi}(\lambda ) } \Big) \frac{\alpha      dx}{x^{1+\alpha}}.
 \eeqnn
 A simple calculation induces that 
$
 V_\mathbf{T}(\lambda) =\Gamma(1-\alpha)c_0^{1+\alpha} |\varPsi(\lambda)|^\alpha = c\cdot |\varPsi(\lambda)|^\alpha .
$
 \qed  
 
 \textit{Proof for Corollary~\ref{MainProp.02}.} 
 When $b=0$, we have   $\mathit{\Psi}(\lambda) = (\lambda /c)^{1/(1+\alpha)} $ and 
 \beqnn
 \mathbf{E}\big[\exp\{-\lambda \mathbf{T}_{X_\zeta}\}\big] =\exp\{ -\zeta \cdot c^{1/(1+\alpha)}  \cdot  \lambda^{\alpha/(1+\alpha)}  \} .
 \eeqnn 
 By Karamata-Tauberian theorem, we have  as $x\to\infty$,
 \beqnn
 \mathbf{P}\{ \mathbf{T}_{X_\zeta}>x \} \sim \frac{\zeta\cdot c^{1/(1+\alpha)}}{\Gamma(1/(1+\alpha))}   \cdot  x^{-\alpha/(1+\alpha)}.
 \eeqnn
  When $b>0$, we have  $\mathit{\Psi}(\lambda)\sim  \lambda /b $ as $\lambda \to 0+$ and hence 
  \beqnn
 1-  \mathbf{E}\big[\exp\{-\lambda \mathbf{T}_{X_\zeta}\}\big] \sim  \frac{\zeta\cdot c}{b^\alpha} \cdot \lambda^{\alpha}. 
  \eeqnn
  Applying  Karamata-Tauberian theorem again, we have  as $x\to\infty$,
  \beqnn
  \mathbf{P}\{ \mathbf{T}_{X_\zeta}>x \} 
  \sim
  \frac{\zeta\cdot c}{b^\alpha\Gamma(1-\alpha)} \cdot x^{-\alpha}.
  \eeqnn
  \qed

 \bibliographystyle{siam}
 
 \bibliography{Reference}

 \end{document}